\documentclass{amsart}
\usepackage{amssymb, amsmath, amsthm,mathtools}
\usepackage{hyperref}

\title{Thomae formula for $2$ Abelian covers of $\mathbb{CP}^1$ } 
\author[Kopeliovich]{by Yaacov Kopeliovich}
\date{23 September 2010 }
\address{Department of Finance University of Connecticut 110 HillSide Rd. Storss CY 06269}
\email{ykopeliovich@yahoo.com}
\newtheorem{thm}{Theorem}[section]
\newtheorem{cor}[thm]{Corollary}
\newtheorem{lem}[thm]{Lemma}

\newtheorem{prop}[thm]{Proposition}
\theoremstyle{definition}
\newtheorem{defn}[thm]{Definition}

\theoremstyle{remark}
\newtheorem{rmk}[thm]{Remark}

\begin{document}
\begin{abstract}
\footnote{This work wouldn't happen if we hadn't attend the special session organized by Broughton Allen. We thank him for this invitation. We also that  Gabino, Mariela Tony Shaska  for interesting discussions and encouragement}
Let $X$ be an Abelian cover  $\mathbb{CP}^{1}$ ramified at $mr$ points, $\lambda_1...\lambda_{mr}.$ we define a class of non positive divisors on $X$ of degree $g-1$ supported in the pre images of the branch points on $X$, such that the  Riemann theta function doesn't vanish on their image in $J(X).$  We obtain a Thomae formula similar to the formulas [BR],[Na],[Z] ,[EG]  and [Ko]. We show that up  to a certain determinant of the non standard periods of $X$, the value of the Riemann theta function at these divisors raised to a high enough power is a polynomial in the branch point of the curve $X.$ Our approach is based on a refinement of Accola's results and  Nakayashiki's approach explained in [Na] for  Abelian covers.
\end{abstract}
\maketitle
\section{Introduction}
We start with a definition: 
\begin{defn}
Let $X_i,g_i$ be a collection of algebraic curves equipped with a ramified mapping $g_i:X_i\mapsto \mathbb{CP}^{1}$
The  fiber product of $X_i,g_i$ is: 
 $X=\left\{p_1,...,p_n|p_i\in X_i, g_i(p_i)=g_j(p_j)\right\}.$
\end{defn}
Now assume that $X_1,,,X_n$ are hyper-elliptic curves, equipped with the usual ramified mapping $f_i:X_i\mapsto \mathbb{CP}^{1}.$ The ramification points of $X_i$ are given through $(\lambda_{i1}...\lambda_{i2m}),1\leq\i\leq n$  and $\lambda_{ij}\neq\lambda_{kl}.$ Let $X$ be the fiber product of $X_1...X_n$ We prove that: 
\begin{thm} \label{three}
Let $v\in Z_2^n$ and  $\beta_{ij}i$ be integer numbers and $0\leq\beta_{ij}\leq 1.$ such that $$\sum_{i=0}^{m}\overline{\beta_{ij}+v_j}=\frac{r}{2} .$$ Then there is a complex number $\alpha$ such that:
\begin{equation}
\theta\left[u\left(\sum_{i=1}^m\sum_{j=1}^n \beta_{ij} \phi^{-1}\left(\lambda_{ij}\right)+K_{z_0}-\sum_{i=1}^{2^n}\infty_i)\right)\right]=
\alpha\sqrt{\det
C}\times{\prod_{\lambda_{ij}\neq \lambda_{kl}}(\lambda_{ij}-\lambda_{kl})}^{q(\beta_{ij},\beta_{kl})+\frac{1}{2}\gamma_{ijkl}}
\end{equation}
and  $q\left(\beta_{ij},\beta_{kl}\right)=\sum_{v\in Z_2^n}\left(\overline{\beta_{ij}+v_j}-\frac{1}{4}\right)\left(\overline{\beta_{kl}+v_l}-\frac{1}{4}\right)$
where,  $\gamma_{ij,kl}=\frac{1}{8}$ if $i=k$ and $\frac{1}{16}$ otherwise.  
\end{thm}

This theorem is a generalization of the work started by [BR],[Na], [EG] and [Ko]. Using methods from String and Quantum field theory Bershadsky and Radul generalized Thomae formula for hyper-elliptic covers to a non singular cyclic covers of the sphere. [Na] gave a more rigorous proof for the formula suggested by [BR] while [EG] modified Nakayashiki's method and treated a special singular case. [Ko] and [Z] obtained generalizations of these formulas to general cyclic covers. In this note we prove the formula we stated above for another type of covers. These covers have an Abelian $2$ group as their automorphism group. They arise naturally as curves above finite fields that enable the construction of  error correcting codes that are associated with algebraic curves [S], and are called Goppa codes. Our goal is to understand whether Thomae type formulas can be carried to a wider class of curves that arise naturally in number theory and other branches like coding theory. \par
Our proof is a generalization of the approach by Nakayashiki in  [Na],[Ko]. We  produce an integrable differential equation that describes the variation of the logarithm of the theta function with respect to the branch points. [Na] constructs certain analytic quantities of the Riemann surfaces locally ( as algebraic expressions supported by local coordinates around the branch points) and compares them to the global expression as derived in [Fa]. Equating the expansions of the global and the local constructions produces the result. We carry this program below. It turns out that the general case of Abelian $2$ covers doesn't differ much from the case considered in  [Na] and [Ko].\par
It has been a major problem in Riemann Surfaces to invert  integrals that define the period of the curve and obtain a far reaching generalizations of trigonometric and elliptic functions for any periods. Riemann carried this program explicitly only to Hyper-elliptic integrals and the general case is wide open. Our firm belief is that finding other natural cases where such inversion is possible will eventually lead to the solution of the generic cases or at least will characterize situations where such solution may occur.\par An inquisitive reader may notice that our methods carry to arbitrary  Abelian ( not just $2$) groups) covers. We prefer to treat the $2$ group Abelian case first for the sake of easiness and clarity of the result. We plan to complete the case of general Abelian groups in the second part of this work. \par
The approach presented here isn't the only one to look for such formulas. In the case of cyclic covers, Hershel Farkas and his collaborators ([EF],[EiF]) reproved Thomae's original result and used classical approach of Riemann to write the branch points as of cyclic covers as ratios of theta functions. [Z] has generalized these results to general cyclic covers. The formulas he obtained are different from the formulas produced by [Ko] and we investigate this gap.
Lastly [KT] used impressive techniques to derive  similar results for general cyclic covers of degree $3.$ Their approach enables them calculate certain constants $\alpha$ explicitly. It is interesting to see if similar approach is valid for the case we consider in this paper.

\section{Non positive divisors on Riemann surface}
Let $X$ be a Riemann surface and assume that $D=\sum d_iz_i$ is a divisor ( not necessarily positive) on it.
\begin{defn}
$H^0\left(X,\mathcal{O}\left(D\right)\right)$ is the collection of functions $f:X\mapsto \mathbb{CP}^1$ on $X$ such that $div(f)\geq D.$ In the [FK] notation this is the space $\mathfrak{R}\left(D\right).$ \\Let $r(D)= dim H^0\left(X,\mathcal{O}\left(D\right)\right).$
\end{defn}
We seek conditions when $\exists E$ a divisor on $X$ such that $E=\sum e_ix_i,$   $e_i\geq 0$ and $D \equiv E.$ Assume that $D$ is not a positive divisor (otherwise you can set $E=D$) Then if $E$ is positive and equivalent to $D$ there exists a non constant function $f$ such that $div(f)=E/D.$ Therefore $f \in H^0\left(X,\mathcal{O}\left(-D\right)\right).$ (That is $f$ is a function such that $div(f) \geq D.)$ Conclude that $r(-D)>0.$ We showed the following: 
]\begin{lem}
Let $D$ be a non positive divisor. Then if there is $E$ a positive divisor such that $E\equiv D$ then $r(-D)>0.$
\end{lem}
Note that because of Jacobi's inversion theorem if $D$ is a divisor such that $r(-D)>0$ there is always a positive divisor $E$ of degree $g(X)-1$ and $D\equiv E.$
Now assume that $deg E=g(X)-1.$ Apply Riemann Roch and conclude that: $r(-D)=i(-D).$ Choose a base point $z_0$ on $X$ and let $u:X\mapsto Jac(X)$ the standard mapping from $X$ into its Jacobian. Let $K_{z_0}$ be the Riemann constant. Then Using Riemann vanishing theorem for theta functions we have the following non vanishing criteria for theta functions:
\begin{lem}
Let $D,$ $degD=g-1$ be a non positive divisor such that $r(-D)=0$ then $\theta\left(u(D)+K_{z_0}\right)\neq 0.$
\end{lem}

\section{Fiber Products of Hyper-elliptic curves}
In this paper we consider a curve $\phi:X\mapsto \mathbb{CP}^1$  that will be a fiber product of hyper-elliptic curves. More precisely let $X_i, i=1,...n$ be  a collection of hyper-elliptic curves given by the equations $y_i^2=\prod_{j=1}^{2m}(x-\lambda_{ji})$ We assume that the polynomials $f_i(x)=\prod_{j=1}^{n}f_i$  satisfy $(f_i,f_j)=1.$ For such polynomials define the curve $X$ in the affine space $A^{n+1}$ such that the point on the curve is $(x,y_1,...y_n)$ and $y_i^2=f_i(x).$ Define the mapping on the sphere as : $(x,y_1,...y_n)\mapsto x.$ Further define the action of the  group $G=Z_2^n$ by sending the unit vector to the automorphism $(x,y_1...y_i...y_n)\mapsto (x,y_1,...-y_i,..y_n).$ Clearly those automorphisms commute and moreover each one is of order $2.$
Let us calculate the genus of this curve ( see [S] for the original calculation)
\begin{thm}
$g(X)= (mn-2)2^{n-1}+1$
\end{thm}
\noindent \textbf{Proof}:\\
We use the Riemann Hurwitz formula to reach the conclusion. By the properties of the mapping it is clear that every $x\in CP^{1}$ has $2^n$ pre-images. The ramifications points of the mapping are the zeros of the polynomials $f_i(x).$ Thus there are $2mn$ points each one has $2^{n-1}$ pre-images with the overall branching of $1$ in each. Thus the genus is: 
\begin{equation}
g(X)=-2^n+1+\frac{2^{n-1}\times 2mn}{2}=2^{n-1}\left(mn-2\right)+1
\end{equation}
$\blacksquare$
For future reference $r$ is the total ramification of the $X.$ Then we have that: $g(X)-1=\frac{r}{2}-2^n$
\section{Holomorphic differentials for $X$}
In the last section we established the genus of $X.$ Let us know find a basis for analytic differentials of $X.$ They are given by the following theorem: 
\begin{thm}
$\frac{x^{l}dx}{y_{i_1}y_{i_2}...y_{i_s}},$ is the basis of differentials of first kind such on $X$ such that: 
\begin{enumerate}
\item $l<ms-2$
\item $1\leq i_1...i_s\leq s$
\end{enumerate} 
\end{thm}
\noindent \textbf{Proof}:\\
Note that the order of each $y_i$ at $\infty$ is $m$ hence the order of  $\prod_{k=1}^s y_{i_k}$ is $ms$. On the other hand as $\infty$ is not  a branch point by construction, hence order of $dx$ is $-2.$, and,   $\frac{x^idx}{\prod_{k=1}^ry_{i_k}}$ is a differential of the first kind as long as $i\leq ms-2.$ 
$\blacksquare$
\section{Dimension of $H^0(X,\mathcal{O}(-D)$}
In this section we will prove a direct dimension formula for certain type of divisors invariant under the action of $Z_2^n.$ Let $\beta_1...\beta_{ij}$ be vectors such $\beta_k=0,1.$ Let $$D=\sum_{l=1}^{2mn}\beta_l\phi^{-1}\left(\lambda_l\right)-\sum_{i=1}^{2^n}\infty_i$$. If $\alpha\in Z_2^n$  and $v=\left(v_1,...v_n\right), i_j=0,1$ $y_\alpha=\prod_{j=1}^ny_i^{v_i}$ 
\begin{thm}
 Let $$\beta_{k\alpha}=\beta_k+Ord_{P_k}\left(y_\alpha\right)(\mod 2).$$  and  assume that: $\sum_{k=1}^{mn}\beta_{k\alpha}=\frac{r}{2^n}-\tau_\alpha$ Then
 $r(-D)=\sum_{k=1}^{2^n}Max(0,\tau_{\alpha})$
\end{thm}
\noindent \textbf{Proof}:\\
By definition of $D$ we have: $deg(D)=\sum_{i=1}^{mn} 2^{n-1}\beta_{i}-2^n=\frac{r}{2}-2^n=g-1.$  $D$ is an invariant divisor under the action of $Z_2^n$ and therefore:  $$H^(X,O(-D))=\sum_\chi D_\chi$$ and $D_\chi$ are the eigen-spaces that correspond to the character $\chi$ hence if $f\in D_\chi.$ $G$ acts trivially on $\frac{f}{y_{\chi}}$ hence, $\frac{f}{y_{\chi}}\in H^{0}\left(-D-div(-y_{\chi})\right).$ is a lift of a function on $\mathbb{CP}^{1}$ with the  divisor $-D-div(y_\chi).$  $y_\chi$ equals to $y_v$ for one of $\alpha$ defined above. For a divisor $D+y_{\chi}=\sum\gamma_i\left(\phi^{-1}\lambda_i\right)-\gamma_{\infty}\sum_{i=1}^{2^n}\infty_i$ we define: 
$$\sigma_{\chi,D}=\sum_{i=1}^{mn}\left[\frac{\gamma_i}{2}\right]\lambda_i+\infty.$$ We show that: 
\begin{lem} 
$H^0(-\sigma_{\chi,D})$ isomorphic to $H^{0}\left(-D-y_{\chi}\right)^G$
\end{lem} 
\noindent \textbf{Proof}:\\
Let $h\in H^0(-\sigma_{\chi,D}), $ if $\hat{h}$ is a lift of this function to $X$ then $\hat{h}\in H^0(\phi^{-1}\left(-\sigma_{\chi,D}\right)).$ Apply the definition of $D+div(y_{\chi})$ and conclude: $$\hat{h}\in H^0(-D-div(y_{\chi}).$$ Now assume that $h_1\in H^0\left(D+div(y_{\chi})\right)$ is a lift from a function on the sphere. by definition of $H^0\left(D+div(y_{\chi})\right)$ we know that $div(h_1)\geq -D-div(y_\chi)$ but it's a  lift from $\mathbb{CP}^{1}$ hence:  $div(h_1)=\sum_{i=1}^n \beta_i\phi^{-1}\lambda_i$ and $\beta_i=0 \mod 2.$ Applying Riemann Roch to $CP^1$ we see that $r(-\sigma_{\chi,D})=Max(\tau_k,0).$ Conclude the result $\blacksquare$\\
Hence we obtain the following corollary similar to the case of cyclic covers: 
\begin{cor} \label{cornonvanishing}
	Let $\beta_i$ be selected such that $\tau_i=0, 0\leq i\leq 1 $ then $$\theta\left[u(\sum_{i=1}^{mn} \beta_i\phi^{-1}(\lambda_{ij}))+ K_{z_0}-u(\sum_{j=1}^{2^n}\infty_j)\right](0,\tau)\neq 0.$$
\end{cor}
\section{Algebraic construction of the Szego Kernel}
Let us recall the definition of the Szego Kernel.
\begin{defn}
For $e\in \mathbb{C}^g.$ if $\theta[e]\neq 0$ define the Szego kernel by the following equation: $$S[e](P,Q)=\frac{\theta[e]\left(u(P-Q)\right)}{\theta[e](0,\tau)E(P,Q)},P,Q\in X.$$
\end{defn}
$E(P,Q)$ is the prime form. $e$ depends only on its image in the Jacobian , $J(X)$. $S[e](P,Q|)$ has the following properties that are well known [F] (p.19,p.123),and [EG2] (Proof of Theorem 4.7):

\begin{itemize}
  \item $S[e](P,Q)$ is a $\left(\frac{1}{2},\frac{1}{2}\right)$ form with a simple pole along the diagonal
  \item $S[e]\left(P,Q\right)$ has divisor $[e-K_{z_0}]$ with respect to variable $Q$
  \item $S[e]\left(P,Q\right)$ has a divisor $[-e-K_{z_0}]$ with respect to variable $P.$
  \item $S[e](P,Q)$ is a unique up to a constant $\left(\frac{1}{2},\frac{1}{2}\right)$ form that satisfies the previous properties
\end{itemize}
Our goal is to generalize the approach of Nakayashiki and construct $S[e](P,Q)$ algebraically for our type of covers. If $v=(v_1,...v_n),v_i=0,1$   we define for each  $\beta=\left(\beta_1...\beta_{mn}\right)$ : 
$$f_{\beta,v}(x)=\prod_{i=1}^{mn}\prod_{j=1}^n\left(x-\lambda_{ij}\right)^{\left\{\frac{\beta_{ij}+v_j }{2}\right\}-\frac{1}{4}}\sqrt{dx}$$
and $\left\{h\right\}$ is the fractional part of $h.$
to give the following expression to the Szego kernel:
\begin{thm} \label{szego}
Let $P=(x_1,y_1),Q=(x_2,y_2)\in X.$ Choose
 $\overline{\beta}=\left(\beta_1...\beta_{mn}\right) \in \mathbb{Z}^{mn}$ be as in \textbf{~\ref{cornonvanishing}}. if  $e = \sum_{i=1}^{mn}u(\beta_i\phi^{-1}(\lambda_i))+K_{z_0}-u(\sum_{i=1}^{2^n}\infty_i)$ Let \begin{equation}\label{szegoeq}
 F_{\widetilde{\beta}}\left(P,Q\right) =\frac{1}{2^n}\frac{\sum_{v\in {Z_2}^n}f_{\beta,v}(x_1)f_{\beta,v}^{-1}(x_2)}{x_2-x_1}\sqrt{dx_1}\sqrt{dx_2}
\end{equation}
Then \begin{equation}
S[e](P,Q)=F_{\widetilde{\beta}}\left(P,Q\right)
\end{equation}
\end{thm}
\noindent \textbf{Proof:\\} We verify that $F_{\widetilde{\beta}}\left(P,Q\right)$ satisfies the properties characterizing $S[e]\left(P,Q\right).$ the RHS of the equation  
 If $P\neq Q$ and $x_1\neq x_2$, than $F_{\widetilde{\beta}}\left(P,Q\right)$ is regular.\\ We check the case when $x_1=x_2$ but $y_{1i}\neq {y_{2i}}.$ i.e., $\exists i, y_{2i}=-{y_{1i}}.$\\ Now, $\frac{\left\{\beta_{k,i}\right\}}{2}=\frac{\beta_i+k-2h_{ik}}{2}.$ Rewrite $F_{\widetilde{\beta}}\left(P,Q\right)$ as:
$$\frac{1}{2^n(x_2-x_1)} \sum_{v\in Z_2^n}\left(\frac{y_{1v}}{y_{2v}}\right)^k\times\left(\frac{x_1-\lambda_i}{x_2-\lambda_i}\right)^{h_{ik}}\prod_{i=1}^{mn}{\left(\frac{x_1-\lambda_i}{x_2-\lambda_i}\right)}^{\frac{\beta_i-\frac{1}{2}}{2}}\sqrt{dx_1}\sqrt{dx_2}$$
	if $y_{2i}=\pm y_{1i}.$ In the limit when $x_1\rightarrow x_2$,  $y_2\rightarrow -y_1.$ and  $\frac{\sum_{v\in Z_2^n}\pm\frac{y_{1i}}{y_{2i}}}{x_2-x_1}\rightarrow 0.$ Therefore $F_{\widetilde{\beta}}(P,Q)$ is regular when $x_1=x_2$ but $y_1\neq y_2.$   
\\
Let us calculate the expansion of $\frac{1}{2^n}\frac{\sum_{v\in Z_2^n}f_{\beta v}(x_1)f_{\beta v}^{-1}(x_2)}{x_2-x_1}\sqrt{dx_1}\sqrt{dx_2}$ as a function of $x_2$ when the expansion is around $x_1.$ Assuming $x_1=x_2$ we get that the leading coefficient is $\frac{1}{2^n}\sum_{v\in Z_2^n} 1=1.$ For the coefficient in $x_2-x_1$ we obtain using the derivative product rule: 
$$
\sum_{v\in Z_2^n}\sum_{i=1}^{2m}\sum_{j=1}^{n}\pm\left\{\frac{\beta_{ij}+v_j}{2}\right\}\times \left(\frac{1}{x_2-\lambda_{ij}}\right)=0
$$
Taking the second derivative according to $x_1.$ to calculate the coefficient of $(x_1-x_2)$ we arrive to the following result:
\begin{prop}
The expansion of  $F_{\widetilde{\beta}}\left(P,Q\right)$ around $P$ a non branch point is:
\begin{equation}
F_{\widetilde{\beta}}\left(P,Q\right) = \frac{\sqrt{dx_1}\sqrt{dx_2}}{x_2-x_1}\left[1+\frac{1}{2}\sum_{i,j=1}^{i,j=2mn}\frac{q(\beta_i,\beta_j)}{(x_2-\lambda_i)(x_2-\lambda_j)}\times(x_1-x_2)^2+...\right]
\end{equation}
where $$q(\beta_i,\beta_j)=\sum_{v_1,v_2\in Z_2^n}\left(\left\{(\beta_{i,v_1}\right\}-\frac{1}{4}\right)\times\left(\left\{(\beta_{j,v_2})\right\}-\frac{1}{4}\right)$$
Where $x_1,x_2$ are the local coordinate around $P,Q$ respectively and $\beta_{i,v}=\overline{\beta}$
\end{prop}
To complete the proof of theorem \textbf(~\ref{szego}) note that $L(P,Q)=F_{\overline{\beta}}(P,Q)-S[e](P,Q)$ are a section of a line bundle $L_{[e-K_{z_0}]}\bigotimes L_{[-e-K_{z_0}]}.$ Because of the expansion of $F_{\overline{\beta}}(P,Q)$ conclude that $L(P,Q)$ is a holomorphic section of the line bundle. But $$H^0\left(L_{[e-K_{z_0}}]\bigotimes L_{[-e-K_{z_0}]}\right)=  H^0\left(L_{[e-K_{z_0}]}\right)\bigotimes H^0\left(L_{[-e-K_{z_0}]}\right)=0$$
and thus $F_{\overline{\beta}}(P,Q)=S[e](P,Q)$ as required.
$\blacksquare$
\begin{rmk}
The above argument is exactly the method adopted in [Na] to prove the claim for the non singular case. See [EG] for a slightly different approach.
\end{rmk}
Based on the the formula given at the beginning of the section [Na] shows the following expansion for $S[e](P,Q)$ in terms of theta functions:
\begin{cor}
The expansion of the Szego kernel can be given in terms of theta functions as follows:
$$S[e](P,Q)=\frac{\sqrt{dx_1}\sqrt{dx_2}}{x_1-x_2}\times\left[1+\sum_{i=1}^g\frac{\partial \log\theta[e]}{\partial z_i}(0)u_i(x_1)(x_1-x_2)+...\right]$$
$u_i(x)$ is the coefficient of $dx_1$ in the expansion of the holomorphic $v_i(x).$
\end{cor}
Comparing the expansions conclude the following result:
\begin{cor}\label{thetadervanish}
$$\frac{\partial \theta[e]}{\partial z_i}(0)=0$$
\end{cor}
The following is obtained by multiplying the expansions:
\begin{lem}
\begin{equation}
S[e](P,Q)S[-e](P,Q)=\frac{dx_1dx_2}{(x_1-x_2)^2}\left[1+\frac{1}{2}\sum_{i,j=1}^m \frac{q\left(\beta_i,\beta_j\right)}{(x_2-\lambda_i)(x_2-\lambda_i)}(x_1-x_2)^2+...\right]
\end{equation}
\end{lem}
\section{Variational Formula for the Period Matrix}
We like to show that:
\begin{thm} \label{variational}
\begin{equation}
\frac{\partial \tau_{jk}(0)}{dt} = \frac{1}{2}\sum_{Q_i\in \phi^{-1}(\lambda_i)}v_k(Q_i)v_j(Q_i)
\end{equation} 
\end{thm} 
\noindent \textbf{Proof:\\}
We define the connection matrix $\sigma$ and $c$ between our set of differentials we defined before and the canonical basis $v_1...v_g(x)$ ($g$- genus of the curve)
Let $P\in X$ and $u$ local coordinate around $P.$ we will need the following definition: 
\begin{defn} 
 $\omega(P,n)$ be a differential of a second kind satisfying the following conditions: 
\begin{enumerate} 
\item $\omega(P,n)$ is holomorphic except the point $P \in X$ where it has a pole of order $n\geq 2.$ At $P$ we have the expansion of the form: 
\begin{equation} 
\omega(P,n)=-\frac{n-1}{u^n}du(1+O(u^n))
\end{equation}
\item $\omega(P,n)$ has zero $a_j$ periods 
\end{enumerate}
\end{defn} 
The differential is defined uniquely if the local coordinate is fixed. For a branch point lying above $\lambda_i$ we will always take the coordinate $\left(z-\lambda_i\right)^{\frac{1}{2}}.$ The following relation holds: 
$$
\int_{b_j}\omega\left(P,n\right)=-\frac{1}{(n-2)!}v_j^{(n-2)}(P)
$$ 
and $v_j^{(n-2)}$ is the coefficient of $u^{n-2}du$ in the expansion of $v_j(x)$ in $u.$
\begin{lem}
if we expand $v_j(x,t)$ as 	
\begin{equation} 
v_l(x,t)=v_l(x)+v_{l1}(x)t+...,
\end{equation} 
then we have that: 
\begin{equation}
v_{l1}\left(x\right)=\sum_{i_1...i_r,j,l}\sigma_{i_1,...i_r,l}\frac{\lambda_i^{\beta-1}}{\prod_{i_1...i_r}y_{i_1}(\lambda_i)...y_{i_r}(\lambda_i)}\omega(Q_i,\alpha+1)
\end{equation}
\end{lem}
\noindent\textbf{Proof:\\}
Let us assume that the branch point we are varying is going to be $\lambda_1.$ Then for each differential $\psi$ from our basis we have the following expansions varying the point $\lambda_1$ and assuming $P=(z,y_1...y_n)$
$\psi_{t}(P)=\psi(P)+\frac{1}{2}\frac{\psi(P)}{z-\lambda_1}t$ if $\psi$ contains $y_1$ and $0$ otherwise. Using the delta function we can write : 
\begin{equation} 
\psi_{t}(P)=\psi(P)+\frac{1}{2}\frac{\psi(P)}{z-\lambda_1}\delta_{\lambda_1}t
\end{equation}
Now we have the relation: 
\begin{equation}
\psi_t(x)=\sum_{l=1}^g\int_{a_l}\psi_t(x)v_l(x,t) 
\end{equation}
We expand both sides to discover that: 
$$ 
\frac{1}{2}\frac{\psi(P)}{z-\lambda_1}\delta_{\lambda_1}=\sum_{l=1}^g\int_{a_l}\psi_0v_{l1}(x)+\int_{a_l}\frac{1}{2}\frac{\psi(P)}{z-\lambda_1}\delta_{\lambda_1}v_{l}(x)
$$
Now define: 
\begin{equation}
\eta(x)=\frac{1}{2}\frac{\psi(P)}{z-\lambda_1}\delta_{\lambda_1}-\sum_{l=1}^g\int_{a_l}\frac{1}{2}\frac{\psi(P)}{z-\lambda_1}\delta_{\lambda_1}v_{l}(x)
\end{equation}
and we have a system of equations for each $\psi:$ 
\begin{equation} 
\sum_{l=1}^gv_{l1}(x)\int_{a_l}\psi=\eta
\end{equation} 
Now $\eta$ has the following properties: 
\begin{enumerate} 
\item $\forall  l \int_{a_l}\eta=0$ 
\item If $\psi=\delta_{\lambda_1}\psi$ then around every point $R_i$ which is the branch point above 
$\lambda_1$ we have that the expansion of $\eta$ is : 
$$
\eta(u)=2\lambda_1^{m}\frac{du}{\psi'(R_i)u^2}+O(1)
$$
\end{enumerate}
Thus we conclude that : 
\begin{equation}
\eta=2\lambda_1^{m}\sum_{R_i\in \phi^{-1}}\frac{du}{\psi'(R_i)u^2}\omega(R_i,2)
\end{equation}
Using the fact that $\int_{a_j}w\eta=c_{j\eta}$ we conclude the result. 
On the other hand we know that $v_{j}=\sum_{\eta}\sigma_{j\eta}\eta.$ Thus comparing the expansion of  $du$ at each $R_i$ ramification point above $\lambda_1$ we discover that: 
$v_{j1}(R_i)=\sum_{\psi}\sigma_\psi'(R_i)2\lambda_1^m\frac{1}{\psi'()R_i)}$ and hence we have that: 
\begin{equation} 
v_{j1}(x)=\sum_{R_i}v_j(R_i)\omega(R_i,2)
\end{equation} 
Integrating versus $b_k$ we get the result.

\section{Algebraic construction for the canonical differential}
We construct the canonical differential algebraically for cyclic covers.
\begin{defn}
The canonical symmetric differential is a $\omega(x,y)$ is a meromorphic  one differential with respect to $x,y \in X$, having a unique pole of second order when $z$ tends to $w$ with a leading expansion coefficient of $1.$ Further for a canonical homology basis $a_i,b_j, 1\leq i,j\leq g$ we have:
$$
\int_{a_i}\omega(x,y)=0
$$
for fixed $y.$
\end{defn}
First we remind the reader of a possible basis for the holomorphic differentials on
$C.$ For $i_1,...i_r=0,1$ let $s_{i_1...i_r}\left(x,y_1...y_n\right)=\prod_{j=1}^ry_{i_j}$ and the point $\left(x,y_{i_1}...y_{i_n}\right)$ is on the curve. Then we showed that: 
\begin{thm}
$\frac{x^l}{s_{i_1,...i_r}}$	 are differentials of first kind on $X$ if $l<mr-2.$
\end{thm}
Let us change the notation slightly: For each $v\in Z_2^n$ and assume that $i_j$ is the collection of coordinates such that $v_{i_j}=1$ and $0$ otherwise. We set: 
\begin{equation} 
s_{v,l}\left(x,y_1...y_n\right)=\frac{x^l}{\prod_{i_1...i_r}y_{i_1}y_{i_2}...y_{i_r}}
\end{equation} 
Let,  $$P_v^v\left(z,w\right)=\sum A_n^v(w)(z-w)^n,$$ such that:
\begin{enumerate}
\item $A_0^{(v}\left(w\right)=\prod_{i_j=i_1,...i_r}^mf_{i_j}(w)$ and $i_j$ is the collection of coordinates such that $v_{i_j}=1$ and $0$ otherwise. 
\item $A_1^{(l)}=\frac{1}{2}\frac{\partial A_0(l)}{\partial w }$
\end{enumerate}
Set:
\begin{equation}\xi_0(x,y)=\frac{dz(x)dz(y)}{\left(z(x)-z(y)\right)^2}
\end{equation}
\begin{equation}\xi_v\left(x,y\right)=\frac{P_v^{(v)}\left(z(x),z(y)\right)dz(x)dz(y)}{s_v(x)s_v(y)\left(z(x)-z(y)\right)^2}
\end{equation}
\begin{equation} \xi\left(x,y\right)=\frac{1}{2^n}\sum_{v\in Z_2^n}\xi_v\left(x,y\right)
\end{equation}
\begin{prop}

\begin{enumerate}
\item $\xi\left(x,y\right)$ is holomorphic outside the diagonal set $\left\{x=y\right\}.$
\item For a non branch point $P\in X$ take $z$ to be a local coordinate around $P.$ Then the expansion in $z(x)$ at $z(y)$ is :
    $$
    \xi\left(x,y\right)=\frac{dz(x)dz(y)}{\left(z(x)-z(y)\right)^2}+
     O \left(\left(z(x)-z(y)\right)^0\right)
    $$
\end{enumerate}
\end{prop}
\noindent \textbf{Proof:\\}
First we check: If $z(P)=z(Q)$ but $P\neq Q$ on $X,$ then $\xi(P,Q)$ is still regular. We will show it for the automorphism that takes $\left(x,y_1...y_n\right)\mapsto \left(x,-y_1,...y_n\right).$  Assume that $P=(p_1,q_1,...q_n)$ and $Q=(p_1,-q_1,...q_n).$  Let us examine the leading term of the expansion of $\xi_l(x,y)$ around $Q.$ By definition of of $\xi_l(x,y) $ $A_0^v(x,y)/s_v(x)s_v(y)=(-1)^{\delta_{q_1}}$ where $\delta_{q_1}=0$ if $q_1$ doesn't appear in the product of $A_0^l(w)$ and $1$ otherwise. 
Hence we will have that: 
\begin{equation} 
\sum_{v}\xi_v(P,Q)_0=\sum_{i_1...i_r}\left(-1\right)^{\delta_{1i_1...i_r}}\frac{dz(P)dz(Q)}{\left(z(P)-z(Q)\right)^2}=0
\end{equation}
and $\delta_{(1i_1,...i_r)}=0$ if $\forall i_r\neq 1$ and $1$ otherwise. ( the assertion follows immediate because the number of vectors that have $0$ or $1$ respectively in their first coordinate is equal)
The same argument is implied if $P$ and $Q$ are ramification points and are different. Once again there will be an automorphism that takes one to the other and their sum will be $0.$
Hence we showed that this form stays regular if  $P\neq Q$ up to order $1.$ Now the Expansion of the function $z(P)$ at the point $z(Q)$ of order $1$ (that is the coefficient of $\frac{1}{z(P)-z(Q)}$ is: 
\begin{equation}
\frac{\partial \frac{ A_1^v(w)(z-w)}{s_v(z)}}{\partial z}+\frac{\partial \frac{A_0^v(w)}{s_v(z)}}{\partial z }=
\frac{A_1^v(w)}{s_v(z)s_{v}(w)}-\sum_{i_1,...i_r}A_0^v(w)\frac{\partial  log s_v(z)}{\partial z}\frac{1}{ s_v(w)}
\end{equation}
In our case we have that $z=w$ and therefore by definition of $A_1(w)$ we rewrite the last expression as: 
\begin{equation} 
\frac{1}{2}\frac{\partial A_0(w)}{\partial w }-A_0\frac{\partial log s_v(w)}{\partial w}=0
\end{equation} 
By the definition of $A_0^l(w)$ and $s_l(w).$
\begin{cor}
$$\omega(x,y)-\xi(x,y)$$ is holomorphic on $X\times X.$
\end{cor}
 As an immediate corollary we infer that we have polynomials:  
$P_{v_1}^{v_2}$ and $v_1\neq v_2$ such that: 
\begin{equation}
\omega(x,y)-\xi(x,y)=\sum_{v_1,v_2}\frac{P_{v_1}^{v_2}\left(z(x),z(y)\right)dz(x)dz(y)}{s_{v_1}(z(x))s_{v_2}(z(y))}
\end{equation}
Where by modifying the definition of $P_{v_1}^{v_2}$ we can exclude the terms $v_1=v_2.$  write: $$P_{v_1}^{v_2}(z,w)=\sum_{j=0}^{d(k)}A_{v_1j}^{v_2}(w)(z-w)^j.$$
Note that $deg_wP_k^{(l)}(z,w)\leq d(N-l).$
Our aim is to show the following proposition:
\begin{prop}\label{temp}
\begin{equation}
\sum_{v_1}\frac{A_{v_{1},2}^{v_1}\left(Q_i\right)}{\frac{\partial f_{v_1}(\lambda_i)}{{\partial w}}}=-\log\det C
\end{equation}
and $C$ is a $g(X)\times g(X)$ period matrix of non normalized form : $$\left(\int_{a_i}z^{j-1}dz/s_{v}(z)\right).$$
\end{prop}
\noindent \textbf{Proof:\\} Let us take a local coordinate $t=\left(z-\lambda_i\right)^{\frac{1}{2}}$ coordinate around a branch point $Q_{ij}\in \phi^{-1}(\lambda_i),.$ The condition $\int_{a_j}\omega(x,y)=0$ is equivalent to the coefficients of the expansion around $Q_i$ in $dt$ of these expressions vanishes.  This is equivalent to
\begin{equation} \omega^{(1)}\left(x\right)=\frac{1}{2^{n-1}}\sum_{{v_1\in V_{\lambda_i}}}\frac{P_{v_1}^{v_1}(z(x),\lambda_i)dz(x)}{s_{v_1}(x)\frac{\partial{s_{v_1}}}{\partial w}(\lambda_i)(z(x)-\lambda_i)^2}+\sum_{\substack{v_2\\ v_1\neq v_2}}\frac{P_{v_2}^{v_1}\left(z(x),\lambda_i\right)dz(x)}{s_{v_2}(x)\frac{\partial s_{v_1}}{\partial w}(\lambda_i)}
\end{equation}
vanishing when we integrate around $a_j.$  Now write the $\sum\omega^l(Q_i)$ Note: 
$$\frac{\partial}{\partial \lambda_i}\frac{dz}{s_{v_1}}=\frac{1}{2}\frac{dz}{s_{v_1}\left(z-\lambda_i\right)}$$
and 
 \begin{equation}
 P_{v_1}^{(v_1)}(z,\lambda_i)=\frac{1}{2}\prod_{j=1}^m\frac{\partial f_{v_1(\lambda_i)}}{\partial w} (z-\lambda_i)+\sum_{j=0}^{}P_{v_1,j+2}^{v_1}(z-\lambda_i)^{j+2}.  
 \end{equation} 
Hence we get that: 
\begin{equation}
\frac{\partial f_{v_1(\lambda_i)}}{\partial w}\frac{\partial }{\partial \lambda_i}\int_{a_h}\frac{dz}{s_{v_1}}+\frac{1}{2}\sum_{j=0}^{deg(v_1)-1}{A^{v_1}_{v_1,j+2}(\lambda_i)}\int_{a_h}(z-\lambda_i)^j\frac{dz}{s_{v_1}}
\end{equation}
$$+\sum_{k=1,k\neq l}\sum_{j=0}^{deg(v_2)-2}{P^{v_2}_{v_1,j}}(\lambda_i)\int_{a_h}(z-\lambda_i)^j\frac{dz}{s_{v_2}}=0
$$
This is a system of linear equations $g(X)\times g(X)$ in variables, $P^{(v_1)}_{v_2,j}.$ The matrix of these equations is the $g(X)\times g(X)$ matrix $B,$ and $$B=\int_{a_h}\left(z-\lambda_i\right)^{j-1}dz/s_v(z)$$ For each $v_1$ define matrices $B_l$ obtain from $B$ by replacing the column $\int_{a_h}\frac{dz}{s_v(z)}, 1 \leq g(X)$ with the column:$\frac{\partial}{\partial \lambda_i}\int_{a_h}\frac{dz}{s_l(z)}.$ Then
by Cramer's rule: $$\frac{P_{v_1,2}^{(v_1)}}{\frac{\partial f_{v_1(\lambda_i)}}{\partial w}}=-\frac{\det B_{v_1}}{\det B}.$$ 
We show that: 
\begin{lem}
\begin{equation} 
\sum_{v_1}\frac{\det B_{v_1}}{\det B }=\frac{\partial}{\partial \lambda_i}\det C 
\end{equation}
\end{lem} 
\textbf{Proof:}\\
Fist observe that for any matrix $H$ whose entries are function of some variable $x$ we have that: 
\begin{equation} 
\frac{\partial}{\partial x} \log(det H )=\frac{\det H_l}{H}
\end{equation} 
$\det  H_l$ is the matrix whose $H$ whose $l$-th column is replaced by: $\frac{\partial h_{il}}{\partial x}.$ Now in our case if we replace the columns of the form $(z-\lambda_i)\frac{dz}{s_{v}}$ and $j>0$ then the determinant of such matrix will be $0.$ Thus the only term that survive are the terms for the column $\frac{dz}{s_{w}}$ and $s_w$ contains $z-\lambda_i$ as a polynomial. Conclude that 
\begin{equation} 
\sum_{v_1}\frac{\det B_{v_1}}{\det B }=\frac{\partial}{\partial \lambda_i}\det B  
\end{equation}
Since we can transform $B$ into the matrix $C$ using elementary operations that don't alter the determinant conclude that: 
\begin{equation}  
\sum_{v_1}\frac{\det B_{v_1}}{\det B }=\frac{\partial}{\partial \lambda_i}\det C 
\end{equation}\\ $\blacksquare$ \\
Before moving on let us calculate the coefficient of the expansion of $\omega(x,y)$ in $dz(x)dz(y)$ as a function of $x$ evaluated at $y.$ We have that for $\xi_{v_1}(x,y)$ this coefficient equals to: 
\begin{equation} 
-\frac{1}{4}\frac{\partial^2}{\partial x^2}\log{f_{v_1}(x)}-\frac{1}{8}\left(\frac{\partial \log f_{v_1}(x)}{\partial x}\right)^2+\frac{P_{v_1,2}^{v_1}}{f_{v_1}}
\end{equation} 
Define the following object we will work closely when showing Thomae:
\begin{defn}
Let $P=(x,y) \in X$ be a non branch point with a local coordinate $z.$ Define: $$G_z(z)=\lim_{y \rightarrow x}\left[\omega(x,y)-\frac{dz(x)dz(y)}{{(z(y)-z(x))}^2}\right]
$$
\end{defn}
Now taking the local coordinate $t={(z-\lambda_i)}^{\frac{1}{2}}$ around the branch point $Q_i$ we have the following corollary:
\begin{cor}\label{atomic}
The coefficient of $dt^2$ in the expansion of $G_z(z)$ near a branch point lying $Q_{ijk},0\leq k\leq 2^{n-1}$ above $\lambda_{ij}$ in $t={\left(z-\lambda_i\right)}^{\frac{1}{2}}$ is:
$$-\sum_{j=1,j\neq i}^m\frac{\gamma_{ijkl}}{\lambda_{ij}-\lambda_{kl}}-\log \det C,$$ where
$$\gamma_{ijkl}=\frac{1}{8}$$ if $i=k$ and $\frac{1}{16}$ if $i\neq k.$
\end{cor}
The next fact that we need is shown by [Fa] corrolary 2.12 or [Na]: 
\begin{prop} 
Let $e$ belong to the Jacobian such that $\theta[e](0,\tau)\neq 0$  then:
\begin{equation}
S[e]\left(x,y\right)=\omega(x,y)+\sum_{i,j=1}^g\frac{\partial^2\log\theta[e](0)}{\partial z_i\partial z_j}v_i(x)v_j(y),
\end{equation}
where $v_i(x),v_j(x)$ are holomorphic differentials on the surface. 
\end{prop}
Hence Passing to the local coordinate $t=\left(z-\lambda_i\right)^{\frac{1}{2}}$ we obtain the following corollary:
\begin{cor}\label{CorYaacov}
The coefficient of $dt^2$ in the Laurent expansion of $G_z(z)$ at a branch point $P_i$ is :
\begin{equation}\label{Formulamain}
\sum_{j=1,j\neq i}^{j=2mn}\frac{q(\beta_i,\beta_j)}{\lambda_i-\lambda_j}-\sum_{r,s=1}^g\frac{\partial^2\log\theta[e](0)}{\partial z_i\partial z_j}{v_r}(P_i){v_s}(P_i),
\end{equation}
${v_r}^{\alpha}(P_i)$ is the coefficient of $dt$ in the expansion of $v_r(x)$ in the local coordinate $t.$
\end{cor}
\section{Thomae formula for Abelian $2^n$ covers of the Sphere}
Now let us show Thomae formula. As in [Na] we write the logarithmic derivative of the theta function on the divisor: $e_{\overline{\beta}}=u(\beta_iP_i)+K_{z_0}-u(\sum_{i=1}^{N}\infty_i)$
$$\frac{\partial\log\theta\left[e_{\overline{\beta}}\right]}{\partial \lambda_i}\left(0,\tau\right)={\frac{d}{dt}}\log\theta_t\left[e_{\overline{\beta}}\right](0)|_{t=0}\left(0,\tau\right)
$$
By the chain rule the last expression is:
$$
{\frac{d}{dt}}\log\theta_t\left[e_{\overline{\beta}}\right](0)|_{t=0}\left(0,\tau\right)=\sum_{1\leq k,r\leq g}\frac{\partial\log\theta\left[e_{\overline{\beta}}\right]}{\partial \tau_{kr}}\left(0\right)\frac{d\tau_{kr}}{dt}
$$
Now use the heat equation to rewrite the last expression as:
$$
\frac{1}{2}\sum_{1\leq k,r \leq g}\frac{1}{\theta\left[e_{\overline{\beta}}\right]\left(0,\tau\right)}\frac{\partial^2\theta\left[e_{\overline{\beta}}\right]}{\partial z_k \partial z_r}\left(0\right)\frac{d\tau_{kr}}{dt}
$$
We showed in (\textbf{~\ref{thetadervanish}})
$$\frac{\partial \theta[e]_{\overline{\beta}}}{\partial z_i}(0)=0$$ The last sum equals:
$$\frac{1}{2}\sum_{1\leq k,r\leq g}\frac{\partial^2\log\theta\left[e_{\overline{\beta}}\right]}{\partial z_k \partial z_r}\left(0\right)\frac{d\tau_{kr}}{dt}.$$
Use theorem \textbf{~\ref{variational}}, and corollaries \textbf{~\ref{CorYaacov}}, \textbf{~\ref{atomic}} to conclude that:
$$\frac{\partial\log\theta\left[e_{\overline{\beta}}\right]}{\partial \lambda_i}\left(0,\tau\right)=\frac{1}{2}\frac{\partial}{\partial \lambda_i}\log \det C+\sum_{j=1,j\neq i}^{mn}\frac{q(\beta_i,\beta_j)}{\lambda_i-\lambda_j}+\frac{1}{2}\sum_{j=1,j\neq i}^{j=mn}
\frac{\gamma_{ij}}{\lambda_i-\lambda_j}$$
Integrate the system of first order differential equations to get the following theorem:
\begin{thm} \label{four}
Let $v\in Z_2^n$ and  $\beta_{ij}$ $i,jj$,  integer numbers and $0\leq\beta_{ij}\leq 1.$ such that $$\sum_{i=0}^{m}\overline{\beta_{ij}+v_j}=\frac{r}{2} .$$ Then there is a complex number $\alpha$ such that:
\begin{equation}
\theta\left[u\left(\sum_{i=1}^m\sum_{j=1}^n \beta_{ij} \phi^{-1}\left(\lambda_{ij}\right)+K_{z_0}-\sum_{i=1}^{2^n}\infty_i)\right)\right]=
\alpha\sqrt{\det
C}\times{\prod_{\lambda_{ij}\neq \lambda_{kl}}(\lambda_{ij}-\lambda_{kl})}^{q(\beta_{ij},\beta_{kl})+\frac{1}{2}\gamma_{ijkl}}
\end{equation}
and  $q\left(\beta_{ij},\beta_{kl}\right)=\sum_{v\in Z_2^n}\left(\overline{\beta_{ij}+v_j}-\frac{1}{4}\right)\left(\overline{\beta_{kl}+v_l}-\frac{1}{4}\right)$
where,  $\gamma_{ij,kl}=\frac{1}{8}$ if $i=k$ and $\frac{1}{16}$ otherwise.  
\end{thm}

\end{document}